\documentclass[journal]{IEEEtran}

\ifCLASSINFOpdf

\else

\fi

\usepackage{epsfig} % for postscript graphics files
\usepackage{amsmath} % assumes amsmath package installed
\usepackage{subfigure}
\usepackage{color}

\usepackage{times} % assumes new font selection scheme installed
\usepackage{graphics} % for pdf, bitmapped graphics files
\usepackage{amssymb}  % assumes amsmath package installed

\numberwithin{equation}{section}
\begin{document}
\title{Supplementary Material to: Realizations of a Special Class of Admittances with Strictly Lower
Complexity than Canonical Forms}

\author{Michael Z. Q. Chen,  ~ Kai  Wang, ~  Zhan Shu, ~ and ~ Chanying Li

\thanks{This work is supported by HKU CRCG 201111159110, NNSFC 61004093, ``973 Program'' 2012CB720200, and the ``Innovative research projects for graduate students in universities of Jiangsu province'' CXLX12\_0200.}
\thanks{M.Z.Q. Chen is with the Department of Mechanical Engineering, The University of Hong Kong, Pokfulam
Road, Hong Kong; mzqchen@hku.hk.}
\thanks{K. Wang is with the School of Automation, Nanjing
University of Science and Technology, Nanjing, P.R. China;
kwang0721@gmail.com.}% <-this % stops a space
\thanks{Z. Shu is with the Electro-Mechanical Engineering Group, Faculty of Engineering and the Environment, University of Southampton, UK; z.shu@soton.ac.uk.}
\thanks{C. Li is with the Academy of Mathematics and Systems Science, Chinese Academy of Sciences,
Beijing, P.R. China; cyli@amss.ac.cn.}
}% <-this % stops a space

\maketitle

\begin{abstract}
%\boldmath
This is supplementary material to ``Realizations of a special class of admittances with strictly lower complexity than canonical forms'' \cite{CWSL13}, which presents the detailed proofs of some results. For more background information, refer to  \cite{Che07}--\cite{Smi02} and references therein.
\end{abstract}


% \fi
\IEEEpeerreviewmaketitle

%%%%%%%%%%%%%%%%%%%%%%%%%%%%%%%%%%%%%%%%%%%%%%%%%%%%%%%%%%%%%%%%%%%%%%%%%%%%%%%%

\section{Introduction}

This report presents the proofs of some results in the paper ``Realizations of a special class of admittances with strictly lower complexity
than canonical forms'' \cite{CWSL13}. It is assumed that the numbering of lemmas, theorems, corollaries, and figures in this document agrees with that in the original paper.

%%%%%%%%%%%%%%%%%%%%%%%%%%%%%%%%%%%%%%%%%%%%%%%%%%%%%%%%%%%%%%%%%%%%%%%%%%%%%%%%

\section{Proof of Lemma 2}

\begin{proof}
{\it Sufficiency.} It suffices to prove that the three inequalities of Lemma~1 hold. For Case~$1$,  $d_0 = 0$ makes the three inequalities be equivalent to $a_0d_1 \geq 0$, $a_0 \geq 0$, and $a_1 - d_1 \geq 0$,
which are obviously satisfied because of $a_1 - d_1 \geq 0$ and the assumption that all the four coefficients be nonnegative.
For Case~$2$, $a_1 = 0$, $d_1 = 0$ make the three inequalities be equivalent to $a_0 - d_0 \geq 0$,
which obviously holds.

{\it Necessity.} Suppose that $Y(s)$ is positive-real with at least one of the four coefficients being zero. If the number of zero coefficients is exactly one, then to ensure the three inequalities of Lemma~1 to hold, it is noted that only $d_0 = 0$ is possible. In this case,
a necessary and sufficient condition for $Y(s)$ to be positive-real is $a_1 - d_1 \geq 0$. If the number of zero coefficients is exactly two, then to guarantee the positive-realness of $Y(s)$ only the cases when (1) $a_0=0$ and $d_0=0$; (2) $a_1=0$ and $d_1=0$; (3) $d_0=0$ and $d_1=0$ are possible. $Y(s)$ is positive-real if and only if $a_1 - d_1 \geq 0$ when $a_0=0$ and $d_0=0$. $Y(s)$ is positive-real if and only if $a_0 - d_0 \geq 0$ when $a_1=0$ and $d_1=0$.  $Y(s)$ must always be positive-real when $d_0=0$ and $d_1=0$, implying that $a_1 - d_1 \geq 0$ must hold. If the number of
zero coefficients is exactly three, then to ensure the positive-realness only the cases when (1) $a_0=0$, $d_0=0$, and $d_1=0$; (2) $a_1=0$, $d_0=0$, and $d_1=0$ are possible. It is obvious that the two cases must always be positive-real, which implies $a_1 - d_1 \geq 0$. If the four coefficients are all zero, then $Y(s)$ is always positive-real, and $a_1-d_1\geq 0$ always holds. Summarizing all the above discussion, the two cases in the theorem are obtained.
\end{proof}

\section{Proof of Theorem 1}

\begin{proof}
Based on Lemma~2, it is obvious
that there are two possible cases. For Case~$1$, since $d_0=0$ and
$a_1 - d_1 \geq 0$, $Y(s)$ can be written as
\begin{equation*}
Y(s) = k\frac{a_0 s^2 + a_1 s + 1}{d_1 s^2 + s},
\end{equation*}
and $R_k = a_0 (a_0 + d_1^2 - a_1d_1)$. If $R_k > 0$, then $Y(s)$ is written as
\begin{equation*}
\begin{split}
Y(s) =& k/s + \left( 1/(ka_0/d_1) \right. \\
&\left.+ 1/\left( ka_0^3s/R_k + ka_0^2(a_1-d_1)/R_k \right) \right)^{-1},
\end{split}
\end{equation*}
which is realizable with one inductor, one capacitor, and at most two resistors. If $R_k < 0$, then $Y(s)$ can be written as
\begin{equation*}
Y(s) = k/s + ka_0/d_1 + (- a_0d_1/(kR_k) - a_0d_1^2 s/(kR_k))^{-1},
\end{equation*}
which is realizable with two inductors and two resistors.
For Case~$2$, since $a_1 = 0$, $d_1 = 0$, and $a_0 - d_0 \geq 0$,
$Y(s)$ can be written as
\begin{equation*}
Y(s) = k\frac{a_0 s^2 + 1}{s(d_0 s^2 + 1)} = \frac{k}{s} + \frac{k(a_0 - d_0)s}{d_0 s^2 + 1}
\end{equation*}
and $R_k = (a_0 - d_0)^2$. Since it is assumed that $R_k \neq 0$, we conclude that $R_k > 0$ because of the positive-realness of $Y(s)$. $Y(s)$ is realizable with two inductors and one capacitor, which is of strictly lower complexity than the canonical network.
\end{proof}

\section{Proof of Lemma 3}

\begin{proof}
Suppose that $Y(s)$ can be realized by the lossless network, then by
\cite{Bah84} the even part of $Y^{-1}(s)$ is equal to zero, that is, $Ev\ Y^{-1}(s) = 0$. It therefore follows that
\begin{equation*}
\begin{split}
Ev\ Y^{-1}(s) &= \frac{1}{2}\left( Y^{-1}(s)+Y^{-1}(-s) \right)  \\
         &= \frac{2s^2((a_0d_1-a_1d_0)s^2 +(d_1 - a_1))}{k(a_0 s^2 + a_1 s + 1)(a_0 s^2 -a_1 s
        +1)}  = 0.
\end{split}
\end{equation*}
Thus, $2s^2((a_0d_1-a_1d_0)s^2 +(d_1 - a_1))=0$ holds for all $s$. Then we
have $a_1-d_1=0$ and $a_0d_1 - a_1 d_0 = 0$, which indicates that
$Y(s)=k/s$.
\end{proof}

\section{Proof of Lemma 4}
\begin{proof}
It is obvious that $Z(s)=Y^{-1}(s)$ has a pole at $s=\infty$ and a
zero at $s=0$. By \cite[Theorem 2]{Ses59}, this lemma can be easily
proven.
\end{proof}

\section{Proof of Lemma 5}

\begin{proof}
For the network in Fig.~3(b), we see that
there must exist poles on the imaginary axis $s=j\omega_0$ with $\omega_0 \neq 0$, which contradicts the
fact that all the coefficients be positive. Since any network in
the form of Fig.~3(a) is the frequency-inverse dual of a network in Fig.~3(b), then it cannot realize this class of admittances, either.
\end{proof}

\section{Proof of Lemma 7}

\begin{proof}
It has been discussed that $Y_1(s)$ can be written as $Y_1(s)=Y^{-1}(s^{-1})=k'(a'_0s^2 + a'_1 s + 1)/\left(s(d'_0 s^2 + d'_1 s  + 1)\right)$ with $a'_0 = 1/d_0$, $a'_1 = d_1/d_0$, $d'_0 = 1/a_0$, and $d'_1 = a_1/a_0$.
Calculating the corresponding $R_{k_1}$, we obtain $R_{k_1} = \left( (a_0 - d_0)^2 - (a_1 - d_1)(a_0 d_1 - a_1
d_0) \right)/(a_0^2d_0^2) = R_k / (a_0^2d_0^2)$.
Thus, this completes the proof.
\end{proof}

\section{Proof of Lemma 8}

{\it Sufficiency.}
Since $R_k = 0$, there must be at least one
common factor between $(a_0 s^2 + a_1 s + 1 )$ and $(d_0 s^2 + d_1 s
+ 1 )$. Therefore, $Y(s)$ must be of the form
\begin{equation*}
Y(s) = k \frac{(A s + 1)(B s + 1)}{s(C s + 1)(B s + 1)} = k \frac{A s + 1}{s(C s + 1)},
\end{equation*}
where $A$, $B$, $C$ $> 0$. Comparing the above equation with
\begin{equation}   \label{eq: admittance}
Y(s) = k \frac{a_0s^2 + a_1 s + 1}{s(d_0 s^2 + d_1 s + 1)},
\end{equation}
where $a_0$, $a_1$, $d_0$, $d_1$ $\geq 0$ and $k > 0$,
we have the following relations: $a_0 = A B$, $a_1 = A + B$, $d_0 = B C$, and $d_1 = B + C$.
Then we obtain $a_0/d_0 = A/C$.
Since $Y(s)$ is positive-real, it is therefore implied that $a_0 - d_0 \geq
0$, which indicates that $A\geq C$. If $A=C$, then $Y(s)$ reduces to
$Y(s)=k/s$, realizable as just an inductor.
Otherwise, $A>C$ leads to
\begin{equation*}
Y(s) = k \frac{A s + 1}{s(C s + 1)} = \frac{k}{s} + \frac{1}{\frac{C s}{k(A-C)} +
     \frac{1}{k(A-C)}},
\end{equation*}
which is realizable with three elements.

{\it Necessity.} By the method of enumeration, one half of network graphs of two-terminal networks with at most three elements are shown in Fig.~S1.
Lemma~4 implies that the network with the network graph shown in Fig.~S1(a) can only realize the network in Fig.~5(a). If we disregard the networks which can always be reduced to those with less elements, then graphs in Fig.~S1(b) and Fig.~S1(c) are immediately eliminated, because
all the elements could only be inductors by Lemma~4. Then based on Lemma~5, we conclude that only the network shown in Fig.~5(b) is possible,
which is equivalent to its frequency-inverse dual as shown in Fig.~5(c). In summary, $Y(s)$ can be realized as the admittance of at least one of the networks shown in Fig.~5. By calculating their admittances, we see that $R_k = 0$.

\begin{figure}[thpb]
      \centering
      \renewcommand\thefigure{S\arabic{figure}}
      \subfigure[]{
      \includegraphics[scale=0.9]{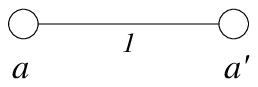}
      \label{subfig: en3(a)}}
      \hspace{0.2cm}
      \subfigure[]{
      \includegraphics[scale=0.9]{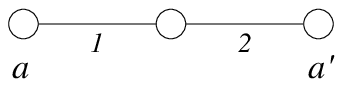}
      \label{subfig: en3(b)}} \\
      \subfigure[]{
      \includegraphics[scale=0.9]{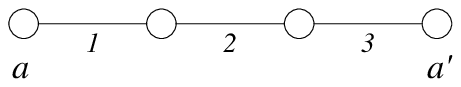}
      \label{subfig: en3(c)}}
      \hspace{0.2cm}
      \subfigure[]{
      \includegraphics[scale=0.9]{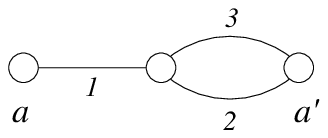}
      \label{subfig: en3(d)}}
      \caption{One half of the network graphs of the two-terminal networks with at most three elements.}
      \label{fig: topological at most three}
   \end{figure}

\section{Proof of Theorem 2}

\begin{proof}
{\it Sufficiency.} For the network shown in Fig.~7(a), its admittance is calculated as
\begin{equation*}
Y(s)=\frac{L_2 C_1 s^2 + R_1 C_1 s + 1}{s(L_1 L_2 C_1 s^2 + R_1
C_1(L_1 + L_2)s + L_1 + L_2)},
\end{equation*}
which can be expressed in the form of \eqref{eq: admittance}, where
$a_0 = L_2C_1 > 0$, $a_1 = R_1C_1 > 0$, $d_0 = L_1L_2C_1/(L_1+L_2) > 0$, $d_1 = R_1C_1 > 0$, and $k = 1/(L_1 + L_2) > 0$. Then it is calculated that $R_k = C_1^2 L_2^4/(L_1+L_2)^2 \neq 0$. So is the admittance of its frequency-inverse dual shown in Fig.~7(b) by Lemma~7, and hence the equivalent networks in Fig.~7(c) and Fig.~7(d). Thus, the sufficiency part of this theorem is proven.

{\it Necessity.} Suppose that $Y(s)$ can be realized with at most four elements whose values are positive and finite. Based on Lemma~8, $R_k \neq 0$ guarantees that we  only need to consider the irreducible four-element network. The method of enumeration is used here for the proof. One half of the  network graphs of the two-terminal four-element networks are listed in Fig.~S2, and other possible graphs are dual with them.

If a four-element network can always be  equivalent to one with less elements, then by Lemma~8 we know $R_k = 0$, contradicting with the assumption. It has been stated in Lemma~4 that there must be a path $\mathcal{P}(a,a')$ and a cut-set $\mathcal{C}(a,a')$ consisting of only inductors for the possible realizations. If the networks whose graphs in Fig.~S2(a), Fig.~S2(b), and Fig.~S2(e) satisfy this property, then they must be equivalent to the ones with at most three elements. Therefore, all these network graphs should be eliminated. For Fig.~S2(c), Edge~$1$ and Edge~$4$ must be inductors by Lemma~4 and  Lemma~8. Since $R_k \neq 0$, $Y(s)$ cannot be written as $Y(s) = k/s$, which means that lossless networks do not need to be considered by Lemma~3. Then there must be exactly one resistor either on Edge~$2$ or on Edge~$3$. If the other element is the inductor, then we obtain the network shown in Fig.~6, which by Lemma~9 means that $R_k =0$. Then, it is eliminated. Therefore, only the network in Fig.~7(a) is possible for the graph in Fig.~S2(c). Using a similar argument, we have that only the network in Fig.~7(d) is possible for graph in Fig.~S2(d). Using the frequency-inverse dual operation, the networks in Fig.~7(b) and Fig.~7(c) are obtained. It is noted that by Lemma~6 the networks in Fig.~7(a) and Fig.~7(b) are equivalent to those in Fig.~7(c) and Fig.~7(d), respectively. So far, all the possible networks have been discovered.
\end{proof}

\begin{figure}[thpb]
      \renewcommand\thefigure{S\arabic{figure}}
      \centering
      \subfigure[]{
      \includegraphics[scale=0.7]{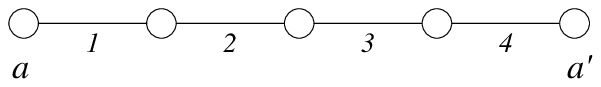}
      \label{subfig: enumerate(a)}}
      \hspace{0.2cm}
      \subfigure[]{
      \includegraphics[scale=0.7]{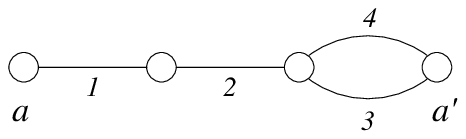}
      \label{subfig: enumerate(b)}} \\
      \subfigure[]{
      \includegraphics[scale=0.7]{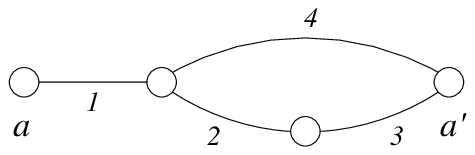}
      \label{subfig: enumerate(c)}}
      \hspace{0.2cm}
      \subfigure[]{
      \includegraphics[scale=0.7]{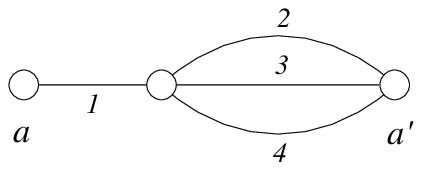}
      \label{subfig: enumerate(d)}}
      \hspace{0.2cm}
      \subfigure[]{
      \includegraphics[scale=0.7]{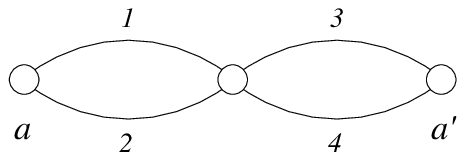}
      \label{subfig: enumerate(e)}}
      \caption{One half of the network graphs of the two-terminal networks with four elements.}
      \label{fig: topological}
   \end{figure}

\section{Proof of Theorem 3}

\textit{Theorem 3:} Consider an admittance $Y(s)$ in the form of \eqref{eq:
admittance} with $a_0$, $a_1$, $d_0$, $d_1$, $k$ $> 0$ and $R_k \neq 0$. Then it can
be realized by the network shown in Fig.~7(a) with
all the elements positive and finite if and only if
\begin{align}
a_0 - d_0 > 0,  \label{eq: four02 condition01}  \\
a_1 - d_1 = 0.  \label{eq: four02 condition02}
\end{align}
Moreover, if the condition is satisfied, then the values of elements are
expressed as
\begin{align}
R_1 &= \frac{a_1 (a_0 - d_0)}{k a_0^2},  \label{eq: four02 R1}  \\
L_1 &= \frac{d_0}{k a_0}, \label{eq: four02 L1}  \\
L_2 &= \frac{a_0 - d_0}{k a_0}, \label{eq: four02 L2} \\
C_1 &= \frac{k a_0^2}{a_0 - d_0}.  \label{eq: four02 C1}
\end{align}

\begin{proof}
{\it Necessity.} The admittance of the network shown in Fig.~7(a) can be calculated as
\begin{equation*}
Y(s) = \frac{L_2C_1 s^2 + R_1C_1 s + 1}{s(L_1L_2C_1s^2 + R_1C_1(L_1
+ L_2)s + L_1 + L_2)}.
\end{equation*}
Suppose that it can realize admittance $Y(s)$ in the form of \eqref{eq: admittance} with $a_0$, $a_1$, $d_0$, $d_1$, $k$ $> 0$, then we obtain
\begin{align}
a_0 &= L_2 C_1,  \label{eq: four02 a0}  \\
a_1 &= R_1 C_1,  \label{eq: four02 a1}  \\
d_0 &= \frac{L_1 L_2 C_1}{L_1 + L_2},  \label{eq: four02 d0}  \\
d_1 &= R_1 C_1,   \label{eq: four02 d1}  \\
k   &= \frac{1}{L_1 + L_2}.  \label{eq: four02 k}
\end{align}
From \eqref{eq: four02 a1} and \eqref{eq: four02 d1}, it is
obvious that \eqref{eq: four02 condition02} holds. From
\eqref{eq: four02 a0} and \eqref{eq: four02 d0}, we obtain
\begin{equation*}
\frac{a_0}{d_0}=\frac{L_1 + L_2}{L_1},
\end{equation*}
from which and together with \eqref{eq: four02 k}, we obtain that $L_1$ can be expressed as \eqref{eq: four02 L1}. Substituting $L_1$ into \eqref{eq: four02 k}, $L_1$ is yielded in the form of \eqref{eq: four02 L2}, from which we indicate \eqref{eq: four02 condition01}. Substituting $L_2$ into \eqref{eq: four02 a0}, we have \eqref{eq: four02 C1}. Finally, from \eqref{eq: four02 a1}, we obtain \eqref{eq: four02 R1}.

{\it Sufficiency.} Consider an admittance $Y(s)$ in the form of \eqref{eq: admittance}, where $a_0$, $a_1$, $d_0$, $d_1$, $k$ $> 0$, $R_k \neq 0$, and \eqref{eq: four02 condition01} and \eqref{eq: four02 condition02} hold. Calculate $R_1$, $L_1$, $L_2$, and $C_1$ by \eqref{eq: four02 R1}--\eqref{eq: four02 C1}. Since \eqref{eq: four02 condition01} holds, it can be verified that all the
values of elements are positive and finite. Since \eqref{eq: four02 condition02} holds,  \eqref{eq:
four02 a0}--\eqref{eq: four02 k} must be satisfied. Therefore, the
admittance of the network in Fig.~7(a) equals
\eqref{eq: admittance}, which implies the sufficiency.
\end{proof}

\section{Proof of Theorem 4}

\begin{proof}
\textit{Sufficiency.} By Lemma~8, the realization of Case~$1$ can obviously be satisfied when $R_k = 0$. It is not difficult to see that the condition of Case~$2$ is equivalent to $a_0 - d_0 > 0$, and either $a_1 - d_1 = 0$ or $a_0d_1 - a_1d_0 = 0$.  When $a_0 - d_0 > 0$ and $a_1 - d_1 = 0$, $Y(s)$ can be realized as in Fig.~7(a) by Theorem~3, which is a four-element network. Through the use of frequency-inverse dual, the case can be proven when $a_0 - d_0 > 0$ and $a_0 d_1 - a_1 d_0 = 0$.

\textit{Necessity.} Suppose that admittance $Y(s)$ in the form of \eqref{eq: admittance} with $a_0$, $a_1$, $d_0$, $d_1$, $k$ $> 0$ can be realized with at most four elements, whose values are positive and finite. We then divide it into two cases: $R_k = 0$ and $R_k \neq 0$. It is known from Lemma~8 that the case when $R_k = 0$ always holds. When $R_k \neq 0$, it is seen from Theorem~2 that $Y(s)$ is the admittance of the network shown in Fig.~7. Together with the realizability conditions derived above, we obtain $a_0 - d_0 > 0$, and either $a_1 - d_1 = 0$ or $a_0 d_1 - a_1 d_0 = 0$, implying Case~2.
\end{proof}

\section{Proof of Theorem 5}

\begin{proof}
\textit{Necessity.} Since $R_k \neq 0$, then $Y(s)$ must be an RL (SD) admittance of McMillan degree three. By  \cite[Theorem 4]{Smi02}, the coefficients must satisfy $R_k < 0$.

\textit{Sufficiency.} Since $R_k:=(a_0 - d_0)^2 - (a_0d_1 - a_1d_0)(a_1 - d_1) < 0$, then we have $a_0d_1 - a_1d_0 > 0$.
In addition, we know $d_0,d_1>0$, then it is known from \cite[Theorem 4]{Smi02} that $Y(s)$ is an RL (SD) admittance of McMillan degree three. From \cite{Smi02,Gui57}, $Y(s)$ can be of the form
\begin{equation*}
Y(s)=k\frac{(As+1)(Cs+1)}{s(Bs+1)(Ds+1)},
\end{equation*}
where $A>B>C>D>0$ and $k>0$. Therefore, it is the admittance of a network with two resistors and an arbitrary number of inductors. Furthermore, $Y(s)$ can be written as
\begin{equation*}
Y(s)=\frac{k}{s} + \frac{k(A-B)(B-C)}{(B-D)(Bs+1)} + \frac{k(A-D)(C-D)}{(B-D)(Ds+1)},
\end{equation*}
which is obviously the admittance of the network in Fig.~8 with the values satisfying $L_1 = 1/k$, $L_2 = B(B-D)/(k(A-B)(B-C))$, $L_3 = D(B-D)/(k(A-D)(C-D))$, $R_1 = (B-D)/(k(A-B)(B-C))$, $R_2 = (B-D)/(k(A-D)(C-D))$,
where $A, C = \left(a_1 \pm \sqrt{a_1^2 - 4a_0}\right)/2$
and $B, D = \left(d_1 \pm \sqrt{d_1^2 - 4d_0}\right)/2$. Since  $a_0 = AC$, $a_1 = A+C$, $d_0=BD$, and $d_1 = B + D$,
then by solving them we obtain the expressions of $A$, $B$, $C$, and $D$ as stated in this theorem.
\end{proof}

\section{Proof of Lemma 10}

\begin{proof}
For Fig.~9(a), the admittance can be calculated as $Y(s)=( C_1L_2L_3s^3+R_1C_1L_3s^2+(L_2+L_3)s+R_1 )/( s( C_1L_1L_2L_3s^3+R_1C_1L_1L_3s^2+(L_1L_2
+L_2L_3+L_1L_3)s+R_1(L_1+L_3) ) )$.
%\begin{equation}  \label{eq: admittance of the network}
%Y(s)=\frac{C_1L_2L_3s^3+R_1C_1L_3s^2+(L_2+L_3)s+R_1}{s\left(C_1L_1L_2L_3s^3+R_1C_1L_1L_3s^2+(L_1L_2
%+L_2L_3+L_1L_3)s+R_1(L_1+L_3)\right)}.
%\end{equation}
Since $R_k\neq 0$, then $Y(s)$ can be realized as in Fig.~9(a), if and only if there exists $T>0$ such that the following equations hold
\begin{align}
a_0 T  &= \frac{C_1L_2L_3}{R_1}, \hspace{0.1cm}
a_0 + a_1 T  = C_1 L_3,         \hspace{0.1cm}
a_1 + T  = \frac{L_2 + L_3}{R_1}, \label{eq: series-parallel 01} \\
&d_0T  = \frac{C_1L_1L_2L_3}{R_1(L_1+L_3)},     \hspace{0.35cm}
d_0 + d_1T  = \frac{C_1L_1L_3}{L_1+L_3}, \label{eq: series-parallel 02} \\
&d_1 + T  = \frac{L_1L_2 + L_2L_3 + L_1L_3}{R_1 (L_1 + L_3)},  \hspace{0.35cm}
k  = \frac{1}{L_1 + L_3}.    \label{eq: series-parallel 03}
\end{align}
Therefore, it suffices to show that $T > 0$ does not exist. After a series of calculations, it is verified that \eqref{eq: series-parallel 01}--\eqref{eq: series-parallel 03} are equivalent to

\begin{equation}
\begin{split}
R_1 &= \frac{(a_0 + a_1 T)(a_0 - d_0)}{ka_0a_1(T^2+a_1T+a_0)} \\
&= \frac{(a_0 - d_0)(a_1d_0T^2 + (a_0^2 + a_1^2d_0)T +a_0a_1d_0)}{ka_0^2a_1(d_1 + T)(T^2 + a_1 T + a_0)}, \label{eq: series-parallel R1}
\end{split}
\end{equation}
\begin{align}
L_1 = &\frac{d_0}{ka_0},  \
L_2 = \frac{(a_0 - d_0)T}{ka_1 (T^2 + a_1T + a_0)},  \
L_3 = \frac{a_0 - d_0}{ka_0},  \label{eq: series-parallel L}
\end{align}
\begin{align}
&C_1 = \frac{ka_0(a_0+a_1T)}{a_0 - d_0} = \frac{ka_0^2(d_0+d_1T)}
{d_0(a_0 - d_0)}.   \label{eq: series-parallel C1}
\end{align}
From \eqref{eq: series-parallel C1}, we obtain $a_0 - d_0 > 0$ and $a_1d_0 - a_0 d_1 = 0$. The second equality of \eqref{eq: series-parallel R1} is equivalent to $a_1(a_0 - d_0)T^2 = 0$,
which obviously cannot hold for any $T>0$. Therefore, there does not exist any $T>0$ such that \eqref{eq: series-parallel 01}--\eqref{eq: series-parallel 03} hold simultaneously.

Using a similar argument, the conclusion of the network in Fig.~9(c) can also be proven. Furthermore, since the networks in Fig.~9(b) and Fig.~9(d) are the frequency-inverse duals of those in Fig.~9(a) and Fig.~9(c), respectively, they cannot be realized, either.
\end{proof}

\section{Proof of Lemma 12}

\textit{Lemma 12:}
Consider a positive-real function $Y(s)$. Then it can be realized as the admittance of the network shown in Fig.~12 with the values of the elements being positive and finite, if and only if $Y(s)$ can be expressed as
\begin{equation}  \label{eq: bridge general admittance}
Y(s)=\frac{\alpha_3 s^3 + \alpha_2 s^2 + \alpha_1 s + 1}{\beta_4 s^4 + \beta_3 s^3 + \beta_2 s^2 + \beta_1 s},
\end{equation}
where $\alpha_1$, $\alpha_2$, $\alpha_2$, $\beta_1$, $\beta_2$, $\beta_3$, $\beta_4$ $> 0$, $W_1$, $W_2$, $W_3$, $W - 2\alpha_2 W_3 > 0$, and the following two equations hold
\begin{align}
W^2 - 4W_1 W_2 W_3 &= 0,  \label{eq: bridge general admittance condition01} \\
\beta_4 + \alpha_1 \beta_3 + \alpha_3\beta_1 - \alpha_2 \beta_2 &= 0. \label{eq: bridge general admittance condition02}
\end{align}
Moreover, if the conditions hold, then the values of the network can be expressed as follows
\begin{align}
R_1 &= \frac{(\alpha_1\alpha_2 - \alpha_3)\beta_1^2}{\alpha_1^2(\alpha_2\beta_1 - \beta_3)}, \label{eq: bridge R1} \\
L_1 &= \frac{(\alpha_1\alpha_2\beta_1 - \alpha_3\beta_1 - \alpha_1\beta_3)\beta_1}{\alpha_1(\alpha_2\beta_1-\beta_3)},  \label{eq: bridge L1} \\
L_2 &= \frac{\alpha_3\beta_1^2}{\alpha_1(\alpha_2\beta_1 - \beta_3)},  \label{eq: bridge L2} \\
L_3 &= \frac{\beta_1\beta_3}{\alpha_2\beta_1-\beta_3}, \label{eq: bridge L3} \\
C_1 &= \frac{\alpha_2\beta_1 - \beta_3}{\beta_1^2}.  \label{eq: bridge C1}
\end{align}

\begin{proof}
\textit{Necessity.} Calculate the admittance of the network in Fig.~12, then we obtain $Y(s) = (C_1L_2(L_1 + L_3)s^3 + R_1 C_1 (L_1 + L_2 + L_3) s^2 + (L_1 + L_3) s + R_1)/(s (C_1L_1L_2L_3 s^3 + R_1C_1L_3(L_1 + L_2) s^2 + (L_1L_2 + L_2L_3 + L_1L_3)s + R_1(L_1 + L_2)))$.
Therefore, $Y(s)$ can be expressed as \eqref{eq: bridge general admittance} with the coefficients satisfying
\begin{align}
&\alpha_3  = \frac{C_1 L_2 (L_1 + L_3)}{R_1}, \hspace{0.35cm}
\alpha_2  = C_1(L_1 + L_2 + L_3),  \label{eq: bridge equations of 01} \\
&\alpha_1  = \frac{L_1 + L_3}{R_1}, \hspace{0.1cm}
\beta_4   = \frac{C_1L_1L_2L_3}{R_1}, \hspace{0.1cm}
\beta_3   = C_1 L_3 (L_1 + L_2),    \label{eq: bridge equations of 02}  \\
&\beta_2   = \frac{L_1L_2 + L_2L_3 + L_1L_3}{R_1},   \hspace{0.35cm}
\beta_1   = L_1 + L_2. \label{eq: bridge equations of 03}
\end{align}
It is obvious that $\alpha_1$, $\alpha_2$, $\alpha_2$, $\beta_1$, $\beta_2$, $\beta_3$, $\beta_4 > 0$. After a series of calculations, the above equations are equivalent to \eqref{eq: bridge R1}--\eqref{eq: bridge C1} and the following two equations
\begin{equation}    \label{eq: bridge equation01}
\begin{split}
\alpha_1\alpha_2^2\beta_1\beta_4 &+ \alpha_3\beta_3\beta_4 + \alpha_3^2\beta_1\beta_3 + \alpha_1\alpha_3\beta_3^2 \\
& = \alpha_2\alpha_3\beta_1\beta_4  + \alpha_1\alpha_2\beta_3\beta_4 + \alpha_1\alpha_2\alpha_3\beta_1\beta_3,
\end{split}
\end{equation}
\begin{equation}   \label{eq: bridge equation02}
\begin{split}
\alpha_1\alpha_2^2&\beta_1\beta_2 + \alpha_3\beta_2\beta_3 + \alpha_3^2\beta_1^2 + \alpha_1\alpha_3\beta_1\beta_3 + \alpha_1^2\beta_3^2 \\
&= \alpha_2\alpha_3\beta_1\beta_2 + \alpha_1\alpha_2\beta_2\beta_3 + \alpha_1\alpha_2\alpha_3\beta_1^2 + \alpha_1^2\alpha_2\beta_1\beta_3.
\end{split}
\end{equation}
Furthermore, \eqref{eq: bridge equation02} is equivalent to
$
\alpha_1\alpha_2\alpha_3\beta_1\beta_3 - \alpha_3^2\beta_1\beta_3 - \alpha_1\alpha_3\beta_3^2 = (\alpha_1\alpha_2^2\beta_1 + \alpha_3\beta_3 - \alpha_2 \alpha_3 \beta_1 - \alpha_1\alpha_2\beta_3)(\alpha_2\beta_2-\alpha_1\beta_3-\alpha_3\beta_1).
$
Since all the values of the elements are positive and finite, then we can calculate that $W_1 = \alpha_1\alpha_2 - \alpha_3 = C_1(L_1 + L_3)^2/R1 > 0$, $W_2 = \alpha_2\beta_1 - \beta_3 = C_1(L_1 + L_2)^2 > 0$, $W_3 = L_1^2/R_1 > 0$, and $W - 2\alpha_2W_3 = 2C_1L_1L_2L_3/R_1 > 0$.
Since $W_1$, $W_2$ $> 0$, then
\begin{equation}  \label{eq: bridge W1W2}
\begin{split}
\alpha_1\alpha_2^2\beta_1 + \alpha_3&\beta_3 - \alpha_2\alpha_3\beta_1 - \alpha_1\alpha_2\beta_3 \\
&= (\alpha_1\alpha_2 - \alpha_3)(\alpha_2\beta_1 - \beta_3) = W_1 W_2 > 0.
\end{split}
\end{equation}
Therefore, it follows from \eqref{eq: bridge equation01} that
\begin{equation}   \label{eq: bridge beta4}
\begin{split}
 \beta_4
=& \frac{\alpha_1\alpha_2\alpha_3\beta_1\beta_3-\alpha_3^2\beta_1\beta_3 - \alpha_1\alpha_3\beta_3^2}{\alpha_1\alpha_2^2\beta_1+\alpha_3\beta_3-
\alpha_2\alpha_3\beta_1-\alpha_1\alpha_2\beta_3}   \\
=& \frac{(\alpha_1\alpha_2^2\beta_1 + \alpha_3\beta_3 - \alpha_2 \alpha_3 \beta_1 - \alpha_1\alpha_2\beta_3)}{\alpha_1\alpha_2^2\beta_1+\alpha_3\beta_3-
\alpha_2\alpha_3\beta_1-\alpha_1\alpha_2\beta_3}  \\
&\times (\alpha_2\beta_2-\alpha_1\beta_3-\alpha_3\beta_1)  \\
=& \alpha_2\beta_2-\alpha_1\beta_3-\alpha_3\beta_1,
\end{split}
\end{equation}
which implies \eqref{eq: bridge general admittance condition02}. Substituting the above equation into $W$,  we obtain $W=2(\alpha_1\alpha_2\beta_1 - \alpha_1\beta_3 - \alpha_3\beta_1)$,
implying
\begin{equation}   \label{eq: bridge W_square - 4W1W2W3}
\begin{split}
4(\alpha_1&\alpha_3\beta_1\beta_3 - \alpha_1^2\alpha_2\beta_1\beta_3 - \alpha_1\alpha_2\alpha_3\beta_1^2 + \alpha_1^2\beta_3^2 + \alpha_3^2\beta_1^2   \\ &+\alpha_1\alpha_2^2\beta_1\beta_2 - \alpha_2\alpha_3\beta_1\beta_2 - \alpha_1\alpha_2\beta_2\beta_3 + \alpha_3\beta_2\beta_3)   \\
=& 4(\alpha_1\alpha_2\beta_1 - \alpha_3\beta_1 - \alpha_1\beta_3)^2 \\
&- 4(\alpha_1 \alpha_2 - \alpha_3)(\alpha_2 \beta_1 - \beta_3)(\alpha_1 \beta_1 - \beta_2)   \\
=& W^2 - 4W_1W_2W_3.
\end{split}
\end{equation}
Now, \eqref{eq: bridge equation02} and \eqref{eq: bridge W_square - 4W1W2W3} yield \eqref{eq: bridge general admittance condition01}.

\textit{Sufficiency.} Let the values of $R_1$, $L_1$, $L_2$, $L_3$, and $C_1$ satisfy \eqref{eq: bridge R1}--\eqref{eq: bridge C1}. $W_1$, $W_2$ $> 0$ indicates that $R_1$, $L_2$, $L_3$, $C_1$ $> 0$. Then substituting $\beta_4$ obtained from \eqref{eq: bridge general admittance condition02} into $W$, we obtain $W-2\alpha_2 W_3 = 2 (\alpha_2\beta_2 - \alpha_1\beta_3 - \alpha_3\beta_1)$. Since $W-2\alpha_2 W_3 > 0$, it follows that $L_1 > 0$. Substituting $\beta_4$ obtained from \eqref{eq: bridge general admittance condition02} into \eqref{eq: bridge general admittance condition01}, we obtain Equation~\eqref{eq: bridge equation02} immediately.  Together with \eqref{eq: bridge W1W2} and \eqref{eq: bridge beta4}, we obtain Equation \eqref{eq: bridge equation01}. It is known in the necessity part that \eqref{eq: bridge equations of 01}--\eqref{eq: bridge equations of 03} must hold.  Now, the sufficiency is proven.
\end{proof}

\section{Proof of Lemma 13}

\textit{Lemma 13:} Consider any positive-real function $Y(s)$ in the form of \eqref{eq: admittance} with $a_0$, $a_1$, $d_0$, $d_1$, $k > 0$, and $R_k \neq 0$. Then it can be realized as the admittance of the network as shown in Fig.~12 with the values of the elements being positive and finite if and only if
\begin{equation}   \label{eq: bridge derived condition}
(a_0d_1 - a_1d_0)(a_1 - d_1) - d_0^2 = 0.
\end{equation}
Moreover, the values of the elements can be expressed as
\begin{align}
R_1 &= \frac{a_1(T^2 + a_1T + a_0)}{k(a_1 + T)^2\left((a_1 - d_1)T + (a_0 - d_0)\right)},  \label{eq: bridge derived R1}   \\
L_1 &= \frac{(a_1 - d_1)T^2 + (a_1^2 - a_1d_1 - d_0)T + a_1(a_0 - d_0)}{k(a_1+T)\left((a_1-d_1)T+(a_0-d_0)\right)},  \label{eq: bridge derived L1} \\
L_2 &= \frac{a_0T}{k(a_1+T)\left((a_1-d_1)T+(a_0-d_0)\right)},  \label{eq: bridge derived L2}  \\
L_3 &= \frac{d_1T + d_0}{k\left((a_1 - d_1)T + (a_0 - d_0)\right)},  \label{eq: bridge derived L3} \\
C_1 &= k\left((a_1-d_1)T + (a_0 - d_0)\right),      \label{eq: bridge derived C1}
\end{align}
where
\begin{equation}  \label{eq: bridge derived T}
T = \sqrt{\frac{a_0d_1-a_1d_0}{a_1-d_1}}.
\end{equation}

\begin{proof}
\textit{Necessity.} Suppose that $Y(s)$ in the form of \eqref{eq: admittance} with all the coefficients positive can be realized as the admittance of the network shown in Fig.~12. It then follows from Lemma~12 that $Y(s)$ can be expressed as \eqref{eq: bridge general admittance} with all the coefficients positive and satisfying the condition of Lemma~12. Since $R_k \neq 0$,  the only way to  express \eqref{eq: admittance} as \eqref{eq: bridge general admittance} where $\alpha_1$, $\alpha_2$, $\alpha_2$, $\beta_1$, $\beta_2$, $\beta_3$, $\beta_4$ $> 0$ is to multiply the numerator and denominator with $(Ts+1)$ where $T>0$. Consequently, it follows that
\begin{equation}   \label{eq: bridge derived coefficients}
\begin{split}
&\alpha_3 = a_0 T,  \  \alpha_2 = a_0 + a_1 T,  \  \alpha_1 = a_1 + T, \\
\beta_4 = &\frac{d_0T}{k},    \
\beta_3 = \frac{d_0 + d_1 T}{k},  \  \beta_2=\frac{d_1 + T}{k}, \ \beta_1 = \frac{1}{k}.
\end{split}
\end{equation}
It is obvious that the condition that $\alpha_1$, $\alpha_2$, $\alpha_2$, $\beta_1$, $\beta_2$, $\beta_3$, $\beta_4$ $> 0$ always holds. Furthermore, other conditions can be presented as follows
\begin{align}
&W_1 = \alpha_1 \alpha_2 - \alpha_3 = a_1T^2 + a_1^2T + a_0a_1 > 0,  \label{eq: bridge derived W1} \\
&W_2 = \alpha_2 \beta_1 - \beta_3 = \frac{(a_1-d_1)T + (a_0-d_0)}{k} > 0,  \label{eq: bridge derived W2} \\
&W_3 = \alpha_1 \beta_1 - \beta_2 = \frac{a_1 - d_1}{k} > 0,     \label{eq: bridge derived W3}            \\
&W - 2\alpha_2W_3 = \frac{(a_1 - d_1)T^2 + (a_0d_1 - a_1d_0)}{k} > 0,  \label{eq: bridge derived W 2alpha2W3} \\
&W^2 - 4W_1W_2W_3 = \frac{\left((a_1-d_1)T^2-(a_0d_1-a_1d_0)\right)^2}{k^2}=0,   \label{eq: bridge derived W_square - 4W1W2W3}
\end{align}
\begin{equation}
\begin{split}
\beta_4 +  \alpha_1 & \beta_3 +  \alpha_3\beta_1 - \alpha_2 \beta_2 \\
&= \frac{-(a_1-d_1)T^2+2d_0T-(a_0d_1-a_1d_0)}{k}=0.  \label{eq: bridge derived equ}
\end{split}
\end{equation}
Then \eqref{eq: bridge derived W3} leads to $a_1 - d_1 > 0$. Then $T$ can be solved from \eqref{eq: bridge derived W_square - 4W1W2W3} as
\eqref{eq: bridge derived T}.
The constraint that $T>0$ yields $a_0d_1 - a_1d_0 > 0$. Substituting the solved $T$ into \eqref{eq: bridge derived equ}, we obtain \eqref{eq: bridge derived condition}.

\textit{Sufficiency.} Suppose that $(a_0d_1 - a_1d_0)(a_1 - d_1) - d_0^2 = 0$ holds, then $d_0 > 0$ and the positive-realness of $Y(s)$ leads to $a_0d_1-a_1d_0>0$ and $a_1-d_1>0$. Then we can let
\begin{equation*}
T = \sqrt{\frac{a_0d_1-a_1d_0}{a_1-d_1}} > 0.
\end{equation*}
Furthermore, multiplying the numerator and denominator of $Y(s)$ with the factor $(Ts+1)$, then we can express $Y(s)$ in the form of \eqref{eq: bridge general admittance} with the coefficients satisfying \eqref{eq: bridge derived coefficients}. It is obvious that $\alpha_1$, $\alpha_2$, $\alpha_2$, $\beta_1$, $\beta_2$, $\beta_3$, $\beta_4 > 0$ and \eqref{eq: bridge derived W1}--\eqref{eq: bridge derived W_square - 4W1W2W3} hold. Since $(a_0d_1 - a_1d_0)(a_1 - d_1) - d_0^2 = 0$, then \eqref{eq: bridge derived equ} is also satisfied. It is concluded that the condition  of Lemma~12 holds, therefore $Y(s)$ can be realized as the admittance of the network as shown in Fig.~12 with the elements positive and finite. The expressions for the values of the elements are presented in \eqref{eq: bridge derived R1}--\eqref{eq: bridge derived C1} derived from \eqref{eq: bridge R1}--\eqref{eq: bridge C1} with the relation \eqref{eq: bridge derived coefficients}.
\end{proof}

\section{Proof of Lemma 14}

\begin{proof}
Calculating the admittance of the network shown in Fig.~13(a) yields $Y_d(s)=(C_1L_2(L_1+L_3)s^3+R_1C_1(L_1+L_3)s^2+(L_1+L_2+L_3)s+R_1)/(s(C_1L_1L_2L_3s^3+R_1C_1
(L_1L_2+L_2L_3+L_1L_3)s^2+L_3(L_1+L_2)s+R_1(L_1+L_2)))$,
which is obviously in the form of \eqref{eq: bridge general admittance} with $\alpha_3 = C_1L_2(L_1+L3)/R_1$, $\alpha_2 = C_1(L_1+L_3)$, $\alpha_1 = (L_1+L_2+L_3)/R_1$, $\beta_4 = C_1L_1L_2L_3/R_1$, $\beta_3 = C_1(L_1L_2+L_2L_3+L_1L_3)$, $\beta_2 = L_3(L_1+L_2)/R_1$, and $\beta_1 = L_1+L_2$.
It can be calculated that the coefficients must satisfy $W-2\alpha_2W_3 = -2C_1L_2(L_1+L_2)(L_1+L_3)/R_1 < 0$. Similarly, for network in Fig.~13(b), its admittance can also be  in the form of \eqref{eq: bridge general admittance} where $\alpha_1$, $\alpha_2$, $\alpha_2$, $\beta_1$, $\beta_2$, $\beta_3$, $\beta_4 > 0$ and satisfy $W - 2\alpha_2W_3 < 0$. Assume that $Y(s)$ in the form of \eqref{eq: admittance} with $a_0$, $a_1$, $d_0$, $d_1$, $k > 0$ and $R_k \neq 0$ can be realized by a network as in Fig.~13(a) or Fig.~13(b). Consequently, $Y(s)$ must be able to be expressed in the form of \eqref{eq: bridge general admittance} with $\alpha_1$, $\alpha_2$, $\alpha_2$, $\beta_1$, $\beta_2$, $\beta_3$, $\beta_4 > 0$, implying that the coefficients satisfy \eqref{eq: bridge derived coefficients}.
In the proof of Lemma~13, it is seen that
\begin{equation*}
W - 2\alpha_2W_3 = \frac{(a_1 - d_1)T^2 + (a_0d_1 - a_1d_0)}{k} \geq 0,
\end{equation*}
which contradicts with the hypothesis. Thus this lemma is proven.
\end{proof}

%%%%%%%%%%%%%%%%%%%%%%%%%%%%%%%%%%%%%%%%%%%%%%%%%%%%%%%%%%%%%%%%%%%%%%%%%%%%%%%%

\section{Conclusion}

In this report, the proofs of some results in the original paper \cite{CWSL13} have been presented.

%%%%%%%%%%%%%%%%%%%%%%%%%%%%%%%%%%%%%%%%%%%%%%%%%%%%%%%%%%%%%%%%%%%%%%%%%%%%%%%%%%%%

\section*{ACKNOWLEDGMENT}
The authors are grateful to the Associate Editor and
the reviewers for their insightful suggestions.


\begin{thebibliography}{99}

\bibitem{CWSL13}
M. Z. Q. Chen, K. Wang, Z. Shu, and C. Li, ``Realizations of a special class of admittances with strictly lower complexity
than canonical forms,'' \textit{IEEE Trans.   Circuits and Systems I: Regular Papers}, vol.~60, no.~9, pp.~2465--2473, 2013.

\bibitem{Che07}
M. Z. Q. Chen, \textit{Passive Network Synthesis of Restricted Complexity}, Ph.D. Thesis, Cambridge Univ. Eng. Dept., U.K., 2007.

\bibitem{CS08}
M.~Z.~Q. Chen and M.~C. Smith, ``Electrical and mechanical passive network synthesis,'' in \textit{Recent Advances in Learning and Control}, V. D. Blondel,  S. P. Boyd, and H. Kimura (Eds.), New York: Springer-Verlag, 2008, LNCIS, vol.~371, pp.~35--50.

\bibitem{Che08}
M.~Z.~Q. Chen, ``A note on PIN polynomials and PRIN rational functions,'' \textit{IEEE Trans. Circuits and Systems II: Express Briefs}, vol.~55, no.~5, pp.~462--463, 2008.

\bibitem{CS09}
M.~Z.~Q. Chen and M.~C. Smith, ``Restricted complexity network realizations for passive mechanical control,'' \textit{IEEE Trans. Automatic Control}, vol.~54, no.~10, pp.~2290--2301, 2009.

\bibitem{CS09(2)}
M.~Z.~Q. Chen and M.~C. Smith, ``A note on tests for positive-real functions,'' \textit{IEEE Trans. Automatic Control}, vol.~54, no.~2, pp.~390--393, 2009.

\bibitem{CPSWS09}
M.~Z.~Q. Chen, C. Papageorgiou, F. Scheibe, F.-C. Wang, and M.~C.
Smith, ``The missing mechanical circuit element,'' \textit{IEEE Circuits Syst. Mag.}, vol.~9, no.~1, pp.~10--26, 2009.


\bibitem{CWZL12}
M.~Z.~Q. Chen, K. Wang, Y. Zou, and J. Lam, ``Realization of a special class of admittances with one damper and one inerter,'' In \textit{Proceedings of the 51st IEEE Conference on Decision and Control}, 2012, pp.~3845--3850.


\bibitem{CWZL13}
M.~Z.~Q. Chen, K. Wang, Y. Zou, and J. Lam, ``Realization of a special class of admittances with one damper and one inerter for mechanical control,'' \textit{IEEE Trans. Automatic Control}, vol.~58, no.~7, pp.~1841--1846, 2013.

\bibitem{CWYLZC13}
M.~Z.~Q. Chen, K. Wang, M. Yin, C. Li, Z. Zuo, and G. Chen,
``Realizability of $n$-port resistive networks with $2n$ terminals,''
in \textit{Proceedings of the 9th Asian Control Conference}, 2013, pp.~1--6.


\bibitem{CWYLZC14}
M.~Z.~Q. Chen, K. Wang, M. Yin, C. Li, Z. Zuo, and G. Chen,
``Synthesis of $n$-port resistive networks containing $2n$ terminals,'' \textit{International Journal of Circuit Theory and Applications}, in press (DOI: 10.1002/cta.1951).

\bibitem{Che14}
M.~Z.~Q. Chen, ``The classical $n$-port resistive synthesis problem,'' in \textit{Workshop on ``Dynamics and Control in Networks''}, Lund University, 2014 (http://www.lccc.lth.se/media/2014/malcolm3.pdf, last accessed on 19/01/2015).



\bibitem{CWZC15}
M.~Z.~Q. Chen, K. Wang, Y. Zou, and G. Chen, ``Realization of three-port spring networks with inerter for effective mechanical control,'' \textit{IEEE Trans. Automatic Control}, in press.


\bibitem{CHD12}
M.~Z.~Q. Chen, Y.~Hu, and B.~Du, ``Suspension performance with one damper and one inerter,''  in \textit{Proceedings of the 24th Chinese Control and Decision Conference}, Taiyuan, China, 2012, pp.~3534--3539.

\bibitem{CHHC14}
M.~Z.~Q. Chen, Y. Hu, L. Huang, and G. Chen, ``Influence of inerter on natural frequencies of vibration systems,'' \textit{Journal of Sound and Vibration}, vol.~333, no.~7, pp.~1874--1887, 2014.

\bibitem{CHLC14}
M.~Z.~Q. Chen, Y. Hu, C. Li, and G. Chen, ``Performance benefits of using inerter in semiactive suspensions,'' \textit{IEEE Trans. Control Systems Technology}, in press (DOI: 10.1109/TCST.2014.2364954).


\bibitem{HCS14}
Y. Hu, M.~Z.~Q. Chen, and Z. Shu, ``Passive vehicle suspensions employing inerters with multiple performance requirements,'' \textit{Journal of Sound and Vibration}, vol.~333, no.~8,
pp.~2212--2225, 2014.


\bibitem{WC12_conf}
K. Wang and M.~Z.~Q. Chen, ``Realization of biquadratic impedances with at most four elements,''  in {\it Proceedings of the 24th Chinese Control and Decision Conference}, 2012, pp.~2900--2905.


\bibitem{WC12}
K. Wang and M.~Z.~Q. Chen, ``Generalized series-parallel RLC synthesis without minimization for biquadratic impedances,''
\textit{IEEE Trans. Circuits and Systems II: Express Briefs}, vol.~59, no.~11, pp.~766--770, 2012.

\bibitem{WCH14}
K. Wang, M.~Z.~Q. Chen, and Y. Hu, ``Synthesis of biquadratic impedances with at most four passive elements,'' \textit{Journal of the Franklin Institute}, vol.~351, no.~3, pp.~1251--1267, 2014.

\bibitem{WC15}
K. Wang and M.~Z.~Q. Chen, ``Minimal realizations of three-port resistive networks,'' \textit{IEEE Trans. Circuits and Systems I: Regular Papers}, in press (10.1109/TCSI.2015.2390560).


\bibitem{Smi02}
M. C. Smith, ``Synthesis of mechanical networks:
the inerter,'' \textit{IEEE Trans. Automatic Control}, vol.~47, no.~10,
pp.~1648--1662, 2002.

\bibitem{Bah84}
H. Baher, \textit{Synthesis of Electrical Networks}. New York: Wiley, 1984.

\bibitem{Ses59}
S. Seshu, ``Minimal realization of the biquadratic minimum
function,'' \textit{IRE Trans. Circuit Theory}, vol.~6, no.~4,
pp.~345--350, 1959.


\bibitem{Gui57} E. A. Guillemin, \textit{Synthesis of Passive
Networks}, John Wiley \& Sons, 1957.



\end{thebibliography}
\end{document}